\documentclass[a4paper,12pt]{article}
\usepackage{amsmath, amssymb, amsthm}
\usepackage{geometry}
\usepackage{graphicx}
\usepackage{xcolor}
\usepackage{hyperref}
\DeclareMathOperator {\tr} {tr} \DeclareMathOperator {\rk} {rk}
\DeclareMathOperator {\Id} {Id}
\newcommand {\const} {\mathrm {const}}
\newcommand {\bbR} {\mathbb R}

\theoremstyle{plain}
\newtheorem{theorem}{\hspace*{\parindent}Theorem}

\newtheorem{remark}{\hspace*{\parindent}Remark}

\let\rom=\textup

\theoremstyle{definition}

\theoremstyle{remark}

\title{Elementary Differential Singularities\\of Three-Dimensional Nijenhuis Operators}
\author{D.Akpan, A.Oshemkov}
\date{}
\begin{document}
\maketitle
\section{Introduction}

Nijenhuis operators are operator fields with zero Nijenhuis torsion. Properties
of Nijenhuis operators, as well as their various applications, have been
actively studied recently. In a recent series of works by A.\,V.~Bolsinov,
A.\,Yu.~Konyaev, and V.\,S.~Matveev, new methods for studying Nijenhuis
operators are developed, and a program for their further research is proposed
\cite{BMK1, BMK2, BMK3, BMK4}. Currently, this program is successfully being
implemented. In particular, some applications of the general theory to various
problems were obtained in the papers~\cite{App1, App2, App3, App4, App5}.

One of the important directions in the development of general theory related to
the Nijenhuis geometry is the study of the singularities of Nijenhuis
operators.

Let us recall the basic definitions (a detailed description of various
properties of Nijenhuis operators can be found in \cite{BMK1}). For an operator
field $L$ (i.e., a tensor field of type $(1,1)$) on a smooth manifold, by its
Nijenhuis torsion one means a tensor of type $(1,2)$ defined by the following
equality:
    $$
    N_L(u,v) = L^2[u,v] + [Lu, Lv] - L[Lu,v] - L[u, Lv],
    $$
where $u,v$ are arbitrary vector fields and $[\,,\,]$ stands for the
commutator. An operator field $L$ is called a Nijenhuis operator if $N_L \equiv
0$.

Let $L$ be an operator field on a smooth manifold $M^n$ and let $\chi(t) = \det
(t\Id - L) = t^n + \sigma_1 t^{n-1} + \ldots + \sigma_n$ be its characteristic
polynomial. Let us define the following mapping given by the invariants
$\sigma_i$ of the operator $L$ regarded as functions on $M^n$:
$$
\Phi_L\colon M^n \to \bbR^n, \qquad \Phi_L(p) = (\sigma_1(p), \sigma_2(p),
\ldots, \sigma_n(p)).
$$
We say that an operator field $L$ has a differential singularity at a point $p
\in M^n$ if $\rk (d \Phi_L) = \rk \bigl(\frac{\partial \sigma_i}{\partial
x^j}\bigr) < n$ at the point $p$. We call such points singular points of rank
$r=\rk(d\Phi_L)$, and operators without singular points are called
differentially nondegenerate.

\begin{theorem}[\cite{BMK1}] \label{theor-can}
Let $L$ be a differentially nondegenerate Nijenhuis operator. Then, at any
local coordinate system $x_1,\dots,x_n$, the following equality holds:
    $$
    L(x) = J(x)^{-1}S_\chi(x) J(x), \quad\text{where}\quad
        S_\chi(x) = \begin{bmatrix}
        -\sigma_1(x)     & 1     & 0      & \ldots & 0 \\
        -\sigma_2(x)     & 0     & 1      & \ldots & 0 \\
        \ldots           & \ldots& \ldots & \ldots & \ldots \\
        -\sigma_{n-1}(x) & 0     & 0      & \ldots & 1 \\
        -\sigma_n(x)     & 0     & 0      & \ldots & 0 \\
    \end{bmatrix},
    $$
and $J(x)=\bigl(\frac{\partial \sigma_i}{\partial x^j}\bigr)$ stands for the
Jacobi matrix of the mapping $\Phi_L$.
\end{theorem}

Theorem \ref{theor-can} gives an explicit formula for the Nijenhuis operator at
any local coordinate system in a neighborhood of the points at which its
invariants $\sigma_1(x),\ldots,\sigma_n(x)$ are functionally independent.
Moreover, any set of functionally independent functions $\sigma_i(x)$ defines
some Nijenhuis operator using the formula from the formulation of
Theorem~\ref{theor-can}. The question concerning the singularities that the
mapping $\Phi_L$ can have, when given by the coefficients of the characteristic
polynomial of some Nijenhuis operator $L$, is still small investigated (this
problem was posed in the paper \cite{BMMT}). In this direction, there is a
paper \cite{AD1} in which two-dimensional Nijenhuis operators are studied for
which the trace is a function without singularities, and the determinant,
restricted to the level line of the trace, has a singularity (for example, a
Morse or a cubic one). Also the result with the Morse singularity was
generalized to arbitrary dimension in the paper \cite{AD2}, where a rather
interesting effect arose: whereas in the two-dimensional case, there were four
Nijenhuis operators with Morse singularity, there are only two of them in the
multidimensional case.

In the present paper, three-dimensional Nijenhuis operators $L$ are considered
for which the corresponding mapping $\Phi_L$ has singularities of some special
type. Namely, we study cases in which $\rk (d\Phi_L)\equiv1$
(Theorem~\ref{theor-rk1}), and also in which $\rk(d\Phi_L)=2$ at some point,
and the two invariants $\sigma_i$ and $\sigma_j$ are functionally independent
in a neighborhood of this point, and the third invariant defines a fold-type
singularity (Theorems~\ref{theor-s1},~~\ref{theor-s2},~~\ref{theor-s3}). As a
result, for singularities of the type under consideration, a complete
description of the corresponding Nijenhuis operators is obtained.

The authors express their gratitude to Andrey Konyaev and Alexey Bolsinov for
valuable advice, as well as to Vladimir Zavyalov, Artemy Sazonov, and Kirill
Dekhnich for support.

\section{Case $\rk(d\Phi_L)\equiv1$.}

In this section, we consider three-dimensional Nijenhuis operators for which
the trace has no singularities and the remaining invariants are functionally
dependent on the trace.

\begin{theorem} \label{theor-rk1}
Let $L$ be a three-dimensional Nijenhuis operator defined on a neighborhood of
the point $p$, and let $\sigma_i$ be the coefficients of its characteristic
polynomial. Suppose that $d\sigma_1\neq 0$ at the point $p$, and the
differentials $d\sigma_2$ and $d\sigma_3$ are proportional to $d\sigma_1$ in a
neighborhood of point~$p$. Then the characteristic polynomial of the operator
$L$ has one of the following three forms\/\rom:

{\rm1. } $t^3-xt^2+(\alpha x-\gamma-\alpha^2)t+\gamma(x-\alpha)$, \
$\alpha,\gamma\in\bbR$,

{\rm2. } $t^3-xt^2+\frac{x^2}3t-\frac{x^3}{27}$,

{\rm3. } $t^3-xt^2+(\frac{x^2}4-\frac{cx}3-\frac{c^2}3)t+\frac
c6(x+\frac{2c}3)^2$, \ $c\in\bbR$,

\noindent in the coordinate system $(x,y,z)$, where $x=-\sigma_1$, and $y,z$
are arbitrary.
\end{theorem}

\begin{proof}
Since $d\sigma_1\ne0$, the trace of the operator $L$ (i.e., the function
$-\sigma_1$) can be taken as the first coordinate: $\sigma_1=-x$. Then the
condition of the theorem exactly means that the remaining invariants depend
only on $x$. For convenience, write $\sigma_2=f(x)$ and $\sigma_3=g(x)$. For
the characteristic polynomial $\chi(L)$ of an arbitrary Nijenhuis operator $L$,
the following formula holds (see~\cite{BMK1}):
$$
L^*(d\chi(t))=t\,d\chi(t)+\chi(t)\,d\tr L.
$$
In the three-dimensional case, in our notation, it can be written out as
follows:
$$
L^*(d(t^3-xt^2+f(x)t+g(x)))=t\,d(t^3-xt^2+f(x)t+g(x))+(t^3-xt^2+f(x)t+g(x))dx,
$$
whence, equating the coefficients of the like powers of~$t$, we obtain
\begin{equation} \label{th2sys1}
    \begin{cases}
        L^*dx = x\,dx-df \\
        L^*df = f\,dx + dg \\
        L^*dg = g\,dx
    \end{cases}
\end{equation}
In particular, from this system it is clear that $L^*dx=a(x)\,dx$, where $a(x)$
is some function, and system \eqref{th2sys1} is rewritten as the following
system of differential equations for the functions $a(x),f(x),g(x)$:
\begin{equation} \label{th2sys2}
    \begin{cases}
        a = x - f' \\
        f'a = f + g' \\
        g'a = g
    \end{cases}
\end{equation}
Differentiating the last equation of system~\eqref{th2sys2} and then
substituting into it expressions for $a$ (from the first equation) and for $g'$
(from the second equation), we obtain
\begin{equation} \label{th2sys3}
(x-f')(x-2f')f''=(f'(x-f')-f)f''.
\end{equation}
If $f''(x) = 0$, then $f(x) = \alpha x + \beta$ and, from
system~\eqref{th2sys2}, we find $a(x) = x - \alpha$, $g(x) = \gamma(x -
\alpha)$, where $\alpha^2 + \beta + \gamma = 0$, which corresponds to the first
case in the formulation of the theorem. If $f''(x)\ne0$, then, from
equation~\eqref{th2sys3}, we obtain
$$
(x - f')(x-3f')+f =0,
$$
which can be rewritten as follows:
$$
(2x-3f')^2=x^2-3f.
$$
The last equation has the form $(F')^2=F$, where $F=x^2-3f$. Its solutions are
$F=0$ and $F=(\frac x2+c)^2$. If $f(x) = \frac{x^2}{3}$, then, from system
\eqref{th2sys2}, we find $a(x) = \frac{x}{3}$ and $g(x) = -\frac{x^3}{27}$.
This is the solution corresponding to case 2 in the formulation of the theorem.
Finally, the last option gives $x^2-3f=(\frac x2+c)^2$, which implies that
$f(x)=\frac{x^2}4-\frac{cx}3-\frac{c^2}3$, and from system \eqref{th2sys2} we
find $a(x)=\frac x2+\frac c3$ and $g(x)=\frac c6(x+\frac{2c}3)^2$, which
corresponds to case 3 in the formulation of the theorem.
     This completes the proof of the theorem.
\end{proof}

\begin{remark}
The characteristic polynomials of Nijenhuis operators in Theorem
\ref{theor-rk1} are factorized as follows:

{\rm1. } $t^3-xt^2+(\alpha
x-\gamma-\alpha^2)t+\gamma(x-\alpha)
=(t-x+\alpha)\Bigl(t-\frac{\alpha-\sqrt{\alpha^2
+4\gamma}}2\Bigr)\Bigl(t-\frac{\alpha+\sqrt{\alpha^2+4\gamma}}2\Bigr)$,

{\rm2. } $t^3-xt^2+\frac{x^2}3t-\frac{x^3}{27}=(t-\frac x3)^3$,

{\rm3. } $t^3-xt^2+(\frac{x^2}4-\frac{cx}3-\frac{c^2}3)t+\frac
c6(x+\frac{2c}3)^2=(t-\frac x2-\frac c3)^2(t+\frac{2c}3)$.

\end{remark}

\section{Case $\rk(d\Phi_L)=2$.}

Let us now consider Nijenhuis operators that have a singular point of rank 2,
and two invariants $\sigma_i$ and $\sigma_j$ are functionally independent, and
the third (after restricting to their joint level line
$\{\sigma_i=0,\,\sigma_j=0\}$) has a Morse singularity at this point. Three
theorems in this section correspond to three possibilities to choose one of the
invariants $\sigma_1,\sigma_2,\sigma_3$ as the one that has Morse
singularity.

    \begin{theorem} \label{theor-s1}
Let $L$ be a three-dimensional Nijenhuis operator defined in a neighborhood of
zero and let $\sigma_i$ be the coefficients of its characteristic polynomial.
Suppose that $d \sigma_2$ and $d \sigma_3$ are linearly independent in the
vicinity of zero, $\sigma_2(0,0,0)=0$ and $\sigma_3(0,0,0) = 0$, and $\sigma_1$
has at zero a Morse singularity on the joint level line of the invariants
$\{\sigma_2=0,\,\sigma_3=0\}$ and $\sigma_1(0,0,0)=c$. Then, for $c=0$, there
are no Nijenhuis operators with the indicated properties, and, if $c\ne0$, then
there exists a regular change of coordinates such that, in the new coordinates,
the matrix of the operator $L$ and its characteristic polynomial have the
following form:
        $$
        L = \begin{bmatrix}
            \mp x^2 - c &  -\frac{x}{2c} &  \frac{x}{2c^2} \\
            \mp 2xy & -\frac{y}{c} & \frac{y}{c^2} + 1 \\
            \mp 2 xz & -\frac{z}{c} & \frac{z}{c^2}
        \end{bmatrix},
        \qquad
        \chi(t) = t^3 + t^2(\pm x^2+\frac yc-\frac z{c^2}+c) + y t + z.
        $$
    \end{theorem}

    \begin{proof}
Since $d\sigma_2$ and $d\sigma_3$ are linearly independent in the vicinity of
zero, it follows that we can take these functions as the second and third
coordinates, respectively; i.e., $\sigma_2 = y$, $\sigma_3 = z$. By Theorem
\ref{theor-can}, we have the following representation for the operator $L$:
                \begin{equation} \label{Lcan-th5}    
{\setlength{\arraycolsep}{2.7pt}
        \! L {=}
        \begin{bmatrix}
            f_x & f_y & f_z \\
            0 & 1 & 0 \\
            0 & 0& 1
        \end{bmatrix}^{-1}      \!
        \begin{bmatrix}
            -f & 1 & 0 \\
            -y & 0 & 1 \\
            -z & 0 & 0
        \end{bmatrix}           \!
        \begin{bmatrix}
            f_x & f_y & f_z \\
            0 & 1 & 0 \\
            0 & 0& 1
        \end{bmatrix}
}
        {=}
        \begin{bmatrix}
            -f + yf_y + zf_z & \frac{f_y(zf_z-f) + yf^2_y + 1}{f_x} &
            \frac{f_y(yf_z-1) + f_z(zf_z-f)}{f_x} \\
            -yf_x & -yf_y & 1 - yf_z \\
            -zf_x & -zf_y & -zf_z
        \end{bmatrix}\!.
        \end{equation}
Since the function $f(x,y,z)$ has a Morse singularity with respect to the
variable $x$, let us use the parametric Morse lemma \cite{AGV} and assume that,
after a regular change of coordinates $(x,y,z)\to(\tilde x(x,y,z),y,z)$, the
function $f$ has the form $f(x,y,z) = \pm x^2 + R(y,z)$ in new coordinates, for
which we keep the same notation. The smoothness of the Nijenhuis operator of
the form \eqref{Lcan-th5} is equivalent to the smoothness of the fractions
$L^1_2$ and $L^1_3$, which have the form
        $$
        \begin{aligned}
            L^1_2&=\frac{f_y(zf_z-f) + yf^2_y + 1}{f_x}
            = -\frac{xR_y}2\pm\frac{R_y(zR_z-R)+yR_y^2+1}{2x},
            \\
            L^1_3&=\frac{f_y(yf_z-1) + f_z(zf_z-f)}{f_x}=
            -\frac{xR_z}2\pm\frac{R_y(yR_z-1)+R_z(zR_z-R)}{2x}.
        \end{aligned}
        $$
From here we obtain the following system of differential equations for the
function $R(y,z)$:
       \begin{equation} \label{sys5}
        \begin{cases}
            R_y(R - zR_z - yR_y) - 1 = 0 \\
            R_z(R - zR_z - yR_y) + R_y = 0
        \end{cases}
        \end{equation}
Subtracting from the first equation multiplied by $R_z$ the second one
multiplied by $R_y$, we obtain
        \begin{equation} \label{Rz}
        R_z=-R^2_y.
       \end{equation}
Let us differentiate the first equation of system \eqref{sys5} with respect to
$y$:
        $$
        R_{yy}(R - zR_z - yR_y) + R_y(-zR_{yz}-yR_{yy})= 0
        $$
Multiplying it by $R_y$ (and taking into account that $R_y(R - zR_z - yR_y) =
1$), we obtain
$$
R_{yy}-R_y^2(zR_{yz}+yR_{yy})=0.
$$
Substituting $R_{yz}=-2R_yR_{yy}$ into this equation (obtained by
differentiating \eqref{Rz} with respect to $y$), we arrive at the equation
        $$
        R_{yy}(1+ 2zR_y^3-yR_y^2)=0.
        $$
The second factor, in some neighborhood of zero, is not equal to zero.
Therefore, we see that $R_{yy}=0$, and also from formula \eqref{Rz} that
$R_{yz}=0$. Thus, $R_y=c_1=\const$. Then from \eqref{Rz} we see that
$R_z=-c_1^2$, i.e., $R=c_1y-c_1^2z+c$. Substituting this expression into the
first equation of system \eqref{sys5}, we obtain
$c_1(c_1y-c_1^2z+c-z(-c_1^2)-yc_1)=1$, whence $c_1=1/c$. Thus, we have proven
that
$$
f(x,y,z)=\pm x^2+\frac yc-\frac z{c^2}+c
$$
and from formula \eqref{Lcan-th5} we obtain the required form of the operator
$L$. This completes the proof of the theorem.
    \end{proof}

\begin{remark}
The characteristic polynomial of the Nijenhuis operators described in
Theorem~\ref{theor-s1}, at the point $p$ under consideration, has the form
$t^2(t-c)$. Therefore, using the splitting theorem from \cite{BMK1}, it can be
shown that, for $c\ne0$, in some coordinate system $u,v,w$, in a neighborhood
of the point $p=(0,0,0)$, the matrices of such operators are reduced to the
form
    $$
    L =
    \begin{bmatrix}
        c \pm u^2 & 0 & 0 \\
        0 & v & 1 \\
        0 & w & 0
    \end{bmatrix}.
    $$
\end{remark}

\begin{theorem} \label{theor-s2}
Let $L$ be a three-dimensional Nijenhuis operator defined in a neighborhood of
zero and let $\sigma_i$ be the coefficients of its characteristic polynomial.
Suppose that $d \sigma_1$ and $d \sigma_3$ are linearly independent in the
vicinity of zero, and also that $\sigma_1(0,0,0)=0$ and $\sigma_3(0,0,0) = 0$.
Let $\sigma_2$ have a Morse singularity on the joint level line of invariants
$\{\sigma_1=0,\,\sigma_3=0\}$. Then in the case $\sigma_2(0,0,0)\ge0$
there are no Nijenhuis operators with the indicated properties, and, if 
$\sigma_2(0,0,0)=-c^2$, where $c\ne0$, then
there is a regular change of coordinates such that, in the new
coordinates, the matrix of the operator $L$ and its characteristic polynomial
have the following form:
    $$
    L =
    \begin{bmatrix}
        -x + c & \pm 2y & \frac{1}{c} \\
        -\frac{y}{2} & -c & 0 \\
        -z & 0 & 0
    \end{bmatrix},
    \quad
    \chi(t) = t^3 + x t^2 + (\pm y^2 + cx + \frac{z}{c} - c^2)t + z.
    $$
\end{theorem}
    \begin{proof}
Since $d\sigma_1$ and $d\sigma_3$ are linearly independent in the vicinity of
zero, it follows that we can take these functions as the first and third
coordinates, respectively; then $\sigma_1 = x$, $\sigma_3 = z$. By Theorem
\ref{theor-can}, we have the following representation for the operator $L$:
        \begin{equation} \label{Lcan-th3}
        L =
        \begin{bmatrix}1&0&0\\f_x&f_y&f_z\\0&0&1\end{bmatrix}^{-1}
        \begin{bmatrix}-x&1&0\\-f&0&1\\-z&0&0\end{bmatrix}
        \begin{bmatrix}1&0&0\\f_x&f_y&f_z\\0&0&1\end{bmatrix}
        =
        \begin{bmatrix}
            -x + f_x & f_y & f_z \\
            \frac{xf_x+zf_z-f-f^2_x}{f_y} & -f_x & \frac{1-f_x f_z}{f_y} \\
            -z & 0 & 0
        \end{bmatrix}.
        \end{equation}
Since the function $f(x,y,z)$ has a Morse (quadratic) singularity with respect
to the variable $y$, we use the parametric Morse lemma \cite{AGV} and assume
that, after a change of coordinates $(x,y,z)\to(x,\tilde y(x,y,z),z)$,
the function $f$ has the form $f(x,y,z) = \pm y^2 + R(x,z)$ in new coordinates,
for which we keep the same notation.

It is clear that the existence of a smooth Nijenhuis operator of the form
\eqref{Lcan-th3} is equivalent to the smoothness of the fractions $L^2_1$,
$L^2_3$, which have the form
$$
        \begin{aligned}
            L^2_1&=\frac{x f_x+z f_z-f- f^2_x}{f_y} =
            \frac{\mp y^2 - R + xR_x+ zR_z - R^2_x}{\pm 2y}
            = -\frac{y}{2} \pm \frac{xR_x+zR_z-R-R^2_x}{2y},
            \\
            L^2_3&=\frac{1-f_x f_z}{f_y} = \pm\frac{1-R_x R_z}{2y}.
        \end{aligned}
        $$
Since the smoothness of a fraction of the form $\frac{F(x,z)}{y}$ in a
neighborhood of the point $(0,0,0)$ is equivalent to $F(x,z)\equiv 0$, we
obtain the following system of differential equations for the function
$R(x,z)$:
        \begin{equation} \label{syst-th4}
        \begin{cases}
            R + R^2_x - xR_x - zR_z = 0\\
            R_x R_z - 1 = 0
        \end{cases}
        \end{equation}
Thus, we are to prove that $R(x,y)$ has a linear form with the corresponding
coefficients. To prove this, let us differentiate the equations with respect to
$x$:
        $$
        \begin{cases}
        2R_xR_{xx} - xR_{xx} - zR_{xz} = 0 \\
        R_{xx}R_z + R_x R_{xz} = 0
        \end{cases}
        \quad\Rightarrow\quad
        \begin{cases}
            2R_xR_{xx} - xR_{xx} - zR_{xz} = 0 \\
            R_{xz} = -R^2_z R_{xx}
        \end{cases}
        $$
        We finally obtain
        $$
        R_{xx}(zR^2_z - x + 2R_x) = 0.
        $$

We have two cases: $R_{xx} = 0$ and $zR^2_z - x + 2R_x = 0$. The second 
case does not give a smooth solution, since it implies $R_x = 0$ 
at the point $(0,0)$, but this contradicts to the equation $R_x R_z = 1$.
Thus, only the equation $R_{xx} = 0$ can give a smooth solution.
If $R_{xx} = 0$, then
        $$
        R(x,z) = x \varphi(z) + \psi(z).
        $$
Therefore
        $$
        R_xR_z = \varphi(z)(x\varphi'(z) + \psi'(z)) = 1.
        $$
This implies
        $$
        \begin{cases}
            \varphi(z) \varphi'(z) = 0 \\
            \varphi(z) \psi'(z) = 1
        \end{cases}
        \quad\Rightarrow\quad
        \begin{cases}
            \varphi(z) = c_1 \neq 0 \\
            \psi(z) = \frac{z}{c_1} + c_2
        \end{cases}
        \quad\Rightarrow\quad
        R(x,z) = c_1x  + \frac{z}{c_1} + c_2,
        $$
where $c_1,c_2$ are some constants.
Since $R(0,0)=-c^2$, we have $c_2=-c^2$, i.e., $R(x,z)=c_1x+\frac{z}{c_1}-c^2$. 
Substituting this into the first equation from \eqref{syst-th4}, we see that
$c_1^2 = c^2$. Thus, we have
        $$
        R(x,z) = cx + \frac{z}{c} - c^2.
        $$
        This completes the proof of the theorem.
    \end{proof}

\begin{remark}
The characteristic polynomial of the Nijenhuis operators described in
Theorem~\ref{theor-s2}, at the point under consideration $p$, has the form
$t(t^2-c^2)$. Therefore, for $c\ne0$, such operators are $gl$-regular and
diagonalizable in some coordinate system in the neighborhood of the point $p$.
\end{remark}

    \begin{theorem} \label{theor-s3}
Let $L$ be a three-dimensional Nijenhuis operator defined in a neighborhood of
zero and let $\sigma_i$ be the coefficients of its characteristic polynomial.
Suppose that $d \sigma_1$ and $d \sigma_2$ are linearly independent in the
vicinity of zero, and also that $\sigma_1(0,0,0)=0$ and $\sigma_2(0,0,0) = 0$.
Let $\sigma_3$ have a Morse singularity on the joint level line of invariants
$\{\sigma_1=0,\,\sigma_2=0\}$ and $\sigma_3(0,0,0)=c$. Then there is a regular
change of coordinates such that, in the new coordinates, the matrix of the
operator $L$ and its characteristic polynomial have the following form:
        $$
        L = \begin{bmatrix}
            - x & 1 & 0 \\
            - y - c^{2/3} & c^{1/3} & \pm 2z \\
            -\frac{z}{2} & 0 & -c^{1/3}
        \end{bmatrix},
        \quad
        \chi(t) = t^3 + xt^2 + yt \pm z^2 + c^{1/3}y - c^{2/3}x + c.
        $$
    \end{theorem}

    \begin{proof}
Since $d \sigma_1$ and $d \sigma_2$ are linearly independent in the vicinity of
zero, it follows that we can take these functions as the first and second
coordinates, respectively; then $\sigma_1 = x$, $\sigma_2 = y$. By Theorem
\ref{theor-can}, we have the following representation for the operator $L$:
        \begin{equation} \label{Lcan-th4}
        L =
        \begin{bmatrix}1&0&0\\0&1&0\\f_x&f_y&f_z\end{bmatrix}^{-1}
        \begin{bmatrix}-x&1&0\\-y&0&1\\-f&0&0\end{bmatrix}
        \begin{bmatrix}1&0&0\\0&1&0\\f_x&f_y&f_z\end{bmatrix}
        =
        \begin{bmatrix}
            - x & 1 & 0 \\
            - y + f_x & f_y & f_z \\
            \frac{xf_x+yf_y -f- f_xf_y}{f_z} & -\frac{f_x + f^2_y}{f_z} & -f_y
        \end{bmatrix}.
        \end{equation}

Since the function $f(x,y,z)$ has a Morse singularity with respect to the
variable $y$, let us use the parametric Morse lemma \cite{AGV} and assume that,
after a regular change of coordinates $(x,y,z)\to(x,y,\tilde z(x,y,z))$, the
function $f$ has the form $f(x,y,z) = \pm z^2 + R(x,y)$ (where $R(0,0) = c$) in
the new coordinates, for which we keep the same notation. The smoothness of the
Nijenhuis operator of the form \eqref{Lcan-th4} is equivalent to the smoothness
of the fractions $L^3_1$ and $L^3_2$, which have the form
$$
        \begin{aligned}
            L^3_1&=\frac{xf_x+yf_y -f- f_xf_y}{f_z} =-\frac{z}{2} \mp
            \frac{R_x R_y - yR_y - xR_x + R}{2z},
            \\
            L^3_2&=-\frac{f_x + f^2_y}{f_z} = \mp\frac{R_x + R^2_y}{2z}.
        \end{aligned}
        $$
Whence we obtain the following system of differential equations for the
function $R(x,y)$:
        $$
        \begin{cases}
            R_x R_y - y R_y - x R_x + R = 0, \\
            R_x + R^2_y = 0.
        \end{cases}
        $$
Let us differentiate the first and second equations with respect to $x$:
        $$
        R_{xx}R_y + R_xR_{xy} - yR_{xy} - xR_{xx} = 0, \quad R_{xx} = -2R_y R_{xy}.
        $$
Then, after substituting $R_{xx}$ into the first equation, we have
        $$
        R_{xy}(R_x - y + 2xR_y - 2R^2_y) = 0,
        $$
that is,
        $$
        R_{xy} = 0 \quad\text{or}\quad
        R_x - y + 2xR_y - 2R^2_y = 0.
        $$
The equation $R_{xy} = 0$ gives the solution $R = \varphi(x) + \psi(y)$; if we
substitute it into the original system, we obtain $\varphi'(x) + (\psi'(y))^2 =
0$. We conclude that the derivatives of $\varphi$ and $\psi$ are constants.
That is, $\varphi'(x) = -c_1$, $(\psi')^2 = c_1$, and hence $\varphi(x) = -c_1
x + c_2$, $\psi(y) = \pm \sqrt{c_1}y + c_3$. Since $R(0,0) = c$, it follows
that $c_2 + c_3 = c$.

Now we substitute the resulting form $R(x,y) = \pm \sqrt{c_1}y - c_1 x + c$
into the first equation of the original system; then we obtain $c_1 = c^{2/3}$.
Thus, we have found the required solution $R(x,y) = c^{1/3}y - c^{2/3}x + c$.

Let us now prove that there are no other solutions, i.e., that $R_x - y + 2xR_y
- 2R^2_y = 0$ does not give a smooth solution. If we substitute $R_x = -R^2_y$,
we obtain
        $$
        3R^2_y - 2xR_y + y = 0
        $$
        $$
        R_y = \frac{x \pm \sqrt{x^2 - 3 y}}{3}
        $$
That is, $R_y$ is not a smooth function.

This completes the proof of the theorem.
    \end{proof}

\begin{remark}
The characteristic polynomial of the Nijenhuis operators described in
Theorem~\ref{theor-s3}, at the point under consideration $p$, has the form
$t^3+c$. Therefore, for $c\ne0$, such operators in a neighborhood of the point
$p$ have one real and two complex conjugate eigenvalues, i.e., are $gl$-regular
and are diagonalizable (in a complex sense) in some coordinate system.
\end{remark}

    \begin{remark}
Theorem 2 in the paper \cite{AD2} for $n = 3$ is a special case of
Theorem~\ref{theor-s3} for $c = 0$.
    \end{remark}

\end{document}